\documentclass[12pt,leqno]{amsart}
\newtheorem{thm}{Theorem}[section]
\newtheorem{defi}{Definition}[section]

\newtheorem{cor}{Corollary}[section]

\newcommand{\be}{\begin{equation}}
\newcommand{\ee}{\end{equation}}
\numberwithin{equation}{section}
\newcommand{\bea}{\begin{eqnarray}}
\newcommand{\eea}{\end{eqnarray}}
\newcommand{\beb}{\begin{eqnarray*}}
\newcommand{\eeb}{\end{eqnarray*}}
\usepackage{amssymb,amsfonts,amsthm,setspace,indentfirst}
\usepackage[dvips]{graphics}
\usepackage{epsfig}

\begin{document}
%%%%%%%%%%%%%%%%%%%%%%%%%%%%%%%%%%%%%%%%%%%%%% Title and address %%%%%%%%%%%%%%%%%%%%%%%%%%%%%%%%%%%%%%%%%%%%%
\title[On Warped Product generalized Roter type Manifold]{\bf{On Warped Product generalized Roter type Manifold}}
\author[Absos Ali Shaikh and Haradhan Kundu]{Absos Ali Shaikh and Haradhan Kundu}
\date{}
\address{\noindent\newline Department of Mathematics,\newline University of 
Burdwan, Golapbag,\newline Burdwan-713104,\newline West Bengal, India}
\email{aask2003@yahoo.co.in, aashaikh@math.buruniv.ac.in}
\email{kundu.haradhan@gmail.com}
\dedicatory{Dedicated to Professor Lajos Tam\'assy on his ninety-first birthday}
%%%%%%%%%%%%%%%%%%%%%%%%%%%%%%%%%%%%%%%%%%%%%%% Abstract and footnote %%%%%%%%%%%%%%%%%%%%%%%%%%%%%%%%%%%%%%%%
%
\begin{abstract}
Generalized Roter type manifold is a generalization of conformally flat manifold as well as Roter type manifold, which gives rise the form of the curvature tensor in terms of algebraic combinations of the fundamental metric tensor and Ricci tensors upto level 2. The object of the present paper is to investigate the characterizations of a warped product manifold to be generalized Roter-type. We also present an example of a warped product manifold which is generalized Roter type but not Roter type, and also an example of a warped product manifold which is Roter type but not conformally flat. These examples ensure the proper existence of such notions.
\end{abstract}
%%%%%%%%%%%%%%%%%%%%%


\maketitle
%
%%%%%%%%%%%%%%%%%%%%%%%%%%%%%%%%%%%%%%%%%%%%%%%%%%%%%%%%%%%%%%%%%%%%%%%%%%%%%%%%%%%%%%%%%%%%%%%%%%%%%%%%%%%%%%%
%                         									    Introduction
%%%%%%%%%%%%%%%%%%%%%%%%%%%%%%%%%%%%%%%%%%%%%%%%%%%%%%%%%%%%%%%%%%%%%%%%%%%%%%%%%%%%%%%%%%%%%%%%%%%%%%%%%%%%%%%
\section{\bf{Introduction}}\label{intro}
%****************************************
The manifold which is locally isometric to an Euclidean manifold is the simplest geometric structure by means of a curvature restriction such that its Riemann-Christoffel curvature tensor $R$ vanishes identically and called a flat manifold. As its proper generalization there arises manifold of constant curvature i.e., the sectional curvatures at each point of the manifold are constant and in this case $R$ is some constant multiple of  the Gaussian curvature tensor $G$ or $g\wedge g$ (for definitions of various symbols used here see Section 2). Conformally flat manifold is a generalization of the manifold of constant curvature, such that $R$ can be expressed as a linear combination of $g\wedge g$ and $g\wedge S$ i.e.,
$$R = J_1 g\wedge g + J_2 g\wedge S,$$
where $J_1$, $J_2$ are some scalars. Especially, for flat manifold $J_1 = J_2 = 0$; for manifold of constant curvature $J_1 = \frac{r}{n(n-1)}$, $J_2 = 0$; and for conformally flat manifold $J_1 = - \frac{r}{2(n-1)(n-2)}$, $J_2 = \frac{1}{n-2}$. Thus we have a way of generalization to find the form of curvature tensor. In this way Roter type manifold (or briefly $RT_n$) \cite{RD03} is a suitable generalization of conformally flat manifold. Similar to conformally flat manifold, in a $RT_n$, the curvature tensor $R$ can be expressed as a linear combination of $g\wedge g$, $g\wedge S$ and $S\wedge S$. Then as a generalization of $RT_n$ in \cite{SDHJK15} Shaikh et. al. introduced the notion of generalized Roter type manifold. A manifold is said to be generalized Roter type (or briefly $GRT_n$) if its curvature tensor is some linear combination of $g\wedge g$, $g\wedge S$, $S\wedge S$ together with $g\wedge S^2$, $S\wedge S^2$ and $S^2\wedge S^2$. We mention that the such decompositions of $R$ were already investigated in \cite{Sawi06} and very recently in \cite{DGJP-TZ13}, \cite{DHJKS13}. We note that the name ``generalized Roter type'' was first used  in \cite{SDHJK15}. For general properties of $GRT_n$ and its proper existence we refer the readers to see \cite{SK} and also references therein.\\
\indent Again the notion of warped product manifold (\cite{BO69}, \cite{Kr57}) is a generalization of product manifold and this notion is important due to its applications in general theory and relativity and cosmology. Various spacetimes are warped product or simply product manifolds, e.g., G\"{o}del spacetime \cite{DHJKS13} is a product manifold and interior black hole spacetime \cite{DHKS}, Robertson-Walker spacetime, generalized Robertson-Walker spacetime are warped products.\\
\indent It is well known that Robertson-Walker spacetime, which is the standard model of cosmology, is a conformally flat warped product space but generalized Robertson-Walker spacetime (\cite{ARS95}, \cite{ARS97}, \cite{EJK96}, \cite{Sa98}, \cite{Sa99}) is a warped product space with $1$-dimensional base and is not conformally flat. In \cite{RD91}, Deszcz studied the conditions for a $4$-dimensional warped product manifolds to be conformally flat (see. Theorem 1, \cite{RD91}) and in \cite{DPS13} Deszcz et. al. showed that certain generalized Robertson-Walker spacetimes are Roter type. In cosmology there arises many non-conformally flat spacetimes which are not Roter type. For example, interior black hole spacetime is non-conformally flat Roter type, and in \cite{DHKS} Deszcz et. al. presented a warped product metric which is not Roter type but of generalized Roter type. We note that in \cite{DPS13} Deszcz et. al. have already studied warped product Roter type manifolds with $1$-dimensional fiber. Motivating by the above studies, in the present paper, we investigate the characterization of a warped product $GRT_n$ manifold.\\
\indent The paper is organized as follows: Section 2 is concerned with preliminaries of such notions. Section 3 deals with warped product manifolds and their different curvature relations. In section 4 we study warped product $GRT_n$ and obtain its characterization (see Theorem \ref{th5.1}). The last section is devoted to the proper existence of such notion with example.
%%%%%%%%%%%%%%%%%%%%%%%%%%%%%%%%%%%%%%%%%%%%%%%%%%%%%%%%%%%%%%%%%%%%%%%%%%%%%%%%%%%%%%%%%%%%%%%%%%%%%%%%%%%%%%%
%                             Preliminaries
%%%%%%%%%%%%%%%%%%%%%%%%%%%%%%%%%%%%%%%%%%%%%%%%%%%%%%%%%%%%%%%%%%%%%%%%%%%%%%%%%%%%%%%%%%%%%%%%%%%%%%%%%%%%%%%
\section{\bf{Preliminaries}}\label{preli}
%****************************************
Let $M$ be an $n (\ge 3)$-dimensional connected semi-Riemannian smooth manifold equipped with a semi-Riemannian metric $g$. We denote by $\nabla$, $R$, $S$, $\kappa$, the Levi-Civita connection, the Riemann-Christoffel curvature tensor, Ricci tensor and scalar curvature of $M$ respectively. The Ricci operator $\mathcal S$ is defined as $g(\mathcal S X, Y) = S(X, Y)$ and the Ricci operator of level 2, $\mathcal S^2$ is defined as $\mathcal S^2 X = \mathcal S(\mathcal S X)$ and its corresponding $(0,2)$ tensor $S^2$, called Ricci tensor of level 2, and is defined as $S^2(X, Y) = S(\mathcal S X, Y)$, where $X, Y \in \chi(M)$, where $\chi(M)$ denotes the Lie algebra of all smooth vector fields on M. In terms of local coordinates the tensor $S^2$ can be expressed as
$$S^2_{ij} = g^{kl}S_{ik}S_{jl}.$$
Similarly we can define the Ricci tensors of level 3 and 4 with corresponding operators as
$$S^3(X,Y) = S(\mathcal S^{2}X,Y), \ \ \ S^2(X,Y) = g(\mathcal S^2 X,Y),$$
$$S^4(X,Y) = S(\mathcal S^{3}X,Y), \ \ \ S^3(X,Y) = g(\mathcal S^3 X,Y).$$
%--------------------------------
Now for $(0,2)$ tensors $A$ and $E$, their Kulkarni-Nomizu product 
(\cite{DG02}, \cite{DGHS98}, \cite{DH03}, \cite{Gl02}) $A\wedge E$ is given by
\bea
(A \wedge E)(X_1,X_2,Y_1,Y_2)&=&A(X_1,Y_2)E(X_2,Y_1) + A(X_2,Y_1)E(X_1,Y_2)\\\nonumber
&&-A(X_1,Y_1)E(X_2,Y_2) - A(X_2,Y_2)E(X_1,Y_1),
\eea
where $X_1, X_2, Y_1, Y_2\in \chi(M)$. 
Throughout the paper we consider $X, Y, X_i, Y_i \in \chi(M)$, $i = 1,2, \cdots $. In particular, we can define $g\wedge g$, $g\wedge S$, $S\wedge S$, $g\wedge S^2$, $S\wedge S^2$ and $S^2\wedge S^2$ etc. The local expression of all such tensors are given by
$$(g\wedge g)_{ijkl} = 2(g_{il}g_{jk} - g_{ik}g_{jl}),$$
$$(g\wedge S)_{ijkl} = g_{il}S_{jk} + S_{il}g_{jk} - g_{ik}S_{jl} - S_{ik}g_{jl},$$
$$(S\wedge S)_{ijkl} = 2(S_{il}S_{jk} - S_{ik}S_{jl}),$$
$$(g\wedge S^2)_{ijkl} = g_{il}S^2_{jk} + S^2_{il}g_{jk} - g_{ik}S^2_{jl} - S^2_{ik}g_{jl},$$
$$(S\wedge S^2)_{ijkl} = S_{il}S^2_{jk} + S^2_{il}S_{jk} - S_{ik}S^2_{jl} - S^2_{ik}S_{jl},$$
$$(S^2\wedge S^2)_{ijkl} = 2(S^2_{il}S^2_{jk} - S^2_{ik}S^2_{jl}).$$
We note that the tensor $\frac{1}{2}(g\wedge g)$ is known as Gaussian curvature tensor and is denoted by $G$.
%------------------------------
A tensor $D$ of type (1,3) on $M$ is said to be generalized curvature tensor (\cite{DG02}, \cite{DGHS98}, \cite{DH03}), if
\beb
&(i)&D(X_1,X_2)X_3+D(X_2,X_1)X_3=0,\\
&(ii)&D(X_1,X_2,X_3,X_4)=D(X_3,X_4,X_1,X_2),\\
&(iii)&D(X_1,X_2)X_3+D(X_2,X_3)X_1+D(X_3,X_1)X_2=0,
\eeb
where $D(X_1,X_2,X_3,X_4)=g(D(X_1,X_2)X_3,X_4)$, for all $X_1,X_2,$ $X_3,X_4$. Here we denote the same symbol $D$ for both generalized curvature tensor of type (1,3) and (0,4). Moreover if $D$ satisfies the second Bianchi identity i.e.,
$$(\nabla_{X_1}D)(X_2,X_3)X_4+(\nabla_{X_2}D)(X_3,X_1)X_4+(\nabla_{X_3}D)(X_1,X_2)X_4=0,$$
then $D$ is called a proper generalized curvature tensor. We note that if $A$ and $B$ are two symmetric $(0,2)$ tensors then $A\wedge B$ is a generalized curvature tensor.\\
%------------------------------
\indent We note that there are various generalized curvature tensors which are linear combination of wedge products of some tensors along with Riemann-Christoffel curvature tensor. One of such more important curvature tensor is the conformal curvature tensor $C$ and is given as
$$C = R -\frac{1}{n-2} g\wedge S + \frac{\kappa}{2(n-1)(n-2)} g\wedge g.$$
We refer the readers to see \cite{SK14} for details about the various curvature tensors and various geometric structures defined through them and their equivalency.
%
%%%%%%%%%%%%%%%%%%%%%%%%%%%%%%%%%%%%%%%%%%%%%%%%%%%%%%%%%%%%%%%%%
%         Roter type and generalized Roter type manifold
%%%%%%%%%%%%%%%%%%%%%%%%%%%%%%%%%%%%%%%%%%%%%%%%%%%%%%%%%%%%%%%%%
\begin{defi}
If the curvature tensor $R$ of a semi-Riemannian manifold $(M^n,g)$ can be expressed as the linear combination of $g\wedge g$, $g\wedge S$ and $S\wedge S$, then it is called Roter type condition on $M$. If a semi-Riemannian manifold satisfies some Roter type condition, then it is called Roter type manifold $($\cite{RD03}$)$ or briefly $RT_n$. Thus on a $RT_n$ we have
\be\label{eqr}
R=N_1 g\wedge g + N_2 g\wedge S + N_3 S\wedge S,
\ee
where $N_1$, $N_2$ and $N_3$ are some scalars on $M$, called the associated scalars of this structure.
\end{defi}
\begin{defi}
If the curvature tensor $R$ of a semi-Riemannian manifold $(M^n,g)$ can be expressed as the linear combination of $g\wedge g$, $g\wedge S$, $S\wedge S$, $g\wedge S^2$, $S\wedge S^2$ and $S^2\wedge S^2$, then this condition is called generalized Roter type condition on $M$. If a semi-Riemannian manifold satisfies some generalized Roter type condition, then it is called generalized Roter type manifold $($\cite{SK}$)$ or briefly $GRT_n$. Thus on a $GRT_n$ we have
\be\label{eqgr}
R=L_1 g\wedge g + L_2 g\wedge S + L_3 S\wedge S + L_4 g\wedge S^2 + L_5 S\wedge S^2 + L_6 S^2\wedge S^2,
\ee
where $L_i$, $1\leqslant i\leqslant 6$ are some scalars on $M$, called the associated scalars of this structure.
\end{defi}
\indent We note that any Roter type manifold is generalized Roter type but not conversely, in general. For details about the geometric properties of generalized Roter type manifold we refer the readers to see \cite{SK}. Throughout this paper by a proper generalized Roter type manifold we mean a $GRT_n$ which is not Roter type, and by a proper Roter type manifold we mean a $RT_n$ which is not conformally flat. A $GRT_n$ or a $RT_n$ is said to be special if one or more of their associated scalars are identically zero.\\
\indent Again by contracting the Roter type and generalized Roter type conditions Shaikh and Kundu \cite{SK} presented various geometric conditions which gives rise to some generalizations of Einstein manifold.
\begin{defi}\cite{Be}
Let $(M^n, g)$ be a semi-Riemannian manifold. If $S$ and $g$ (resp., $S^2$, $S$ and $g$; $S^3$, $S^2$, $S$ and $g$; $S^4$, $S^3$, $S^2$, $S$ and $g$) are linearly dependent, then it is called $Ein(1)$ (resp., $Ein(2)$; $Ein(3)$; $Ein(4)$) condition. A semi-Riemannian manifold satisfying $Ein(i)$ is called $Ein(i)$ manifold, for all $i =1,2,3,4$. The $Ein(1)$ condition is the Einstein metric condition and a manifold satisfying Einstein metric condition is called Einstein manifold and in this case we have $S = \frac{\kappa}{n} g$.
\end{defi}
We note that every $Ein(i)$ manifold is $Ein(i+1)$ for $i=1,2,3$ but not conversely. By taking contraction of the condition of a manifold of constant curvature, Roter type and generalized Roter type, we get the Einstein, $Ein(2)$ and $Ein(4)$ condition respectively. It is well known that every manifold of constant curvature is always Einstein. But a $RT_n$ is $Ein(2)$ except $N_1 = -\frac{\kappa}{2 \left(n^2-3 n+2\right)}$, $N_2 = \frac{1}{n-2}$, $N_3 = 0$; and a $GRT_n$ is $Ein(4)$ except $L_1 = \frac{1}{2} \left(\frac{L_4\left(\kappa^2-\kappa^{(2)}\right)}{n-1}-\frac{\kappa}{n^2-3 n+2}\right)$, $L_2 = \frac{1}{n-2}-L_4 \kappa$, $L_3 = \frac{1}{2} L_4 (n-2)$, $L_5 = 0$, $L_6 = 0$, where $\kappa^{(2)} = tr(S^2)$. We note that the first one gives the conformally flatness.
% 
%%%%%%%%%%%%%%%%%%%%%%%%%%%%%%%%%%%%%%%%%%%%%%%%%%%%%%%%%%%%%%%%%%%%%%%%%%%%%%%%%%%%%%%%%%%%%%%%%%%%%%%%%%%%%%%
%                            Warped Product Manifold
%%%%%%%%%%%%%%%%%%%%%%%%%%%%%%%%%%%%%%%%%%%%%%%%%%%%%%%%%%%%%%%%%%%%%%%%%%%%%%%%%%%%%%%%%%%%%%%%%%%%%%%%%%%%%%%
\section{\bf{Warped Product Manifold}}\label{warp}
%****************************************
The study of warped product manifolds was initiated by Kru$\breve{\mbox{c}}$kovi$\breve{\mbox{c}}$ \cite{Kr57}. Again while  constructing a large class of complete manifolds of negative curvature Bishop and O'Neill \cite{BO69} obtained the notion of the warped product manifolds. The notion of warped product is a generalization of the product of semi-Riemannian manifolds. Let $(\overline M, \overline g)$ and $(\widetilde M, \widetilde g)$ be two semi-Riemannian manifolds of dimension $p$ and $(n-p)$ respectively ($1\leq p < n$), and $f$ be a positive smooth  function on $\overline M$. Let $\overline M$ and $\widetilde M$ be covered with coordinate charts $\left(U; x^1,x^2, ..., x^p\right)$ and $\left(V; y^1, y^2,..., y^{n-p}\right)$ respectively. Then the warped product $M= \overline M\times_{f}\widetilde M$ is the product manifold $\overline M\times \widetilde M$  of dimension $n$ furnished with the metric 
$$g=\pi^*(\overline g) + (f\circ\pi) \sigma^* (\widetilde g),$$
where $\pi:M\rightarrow\overline M$ and $\sigma:M\rightarrow\widetilde M$ are natural projections such that $M = \overline M\times \widetilde M$ is covered with the coordinate charts $\left(U \times V; x^1,x^2, ..., x^p,x^{p+1}=y^1,x^{p+2}=y^2, ...,x^{n}=y^{n-p}\right)$. Then the local components of the metric $g$ with respect to this coordinate chart are given by:
\begin{eqnarray}\label{eq4.1}
g_{ij}=\left\{\begin{array}{lll}
&\overline g_{ij}&\ \ \ \ \mbox{for} \ i = a \ \mbox{and} \ j = b,\\
&f \widetilde g_{ij}&\ \ \ \ \mbox{for $i = \alpha$ and $j = \beta$,}\\
&0&\ \ \ \ \mbox{otherwise.}\\
\end{array}\right.
\end{eqnarray}
Here $a,b \in \left\{1,2,...,p\right\}$ and $\alpha, \beta \in \left\{p+1,p+2,...,n\right\}$. We note that throughout the paper we consider $a,b,c,...\in \{1,2, ..., p\}$ and $\alpha,\beta,\gamma,...\in \{p+1,p+2,...,n\}$ and $i,j,k,...\in \{1,2,...,n\}$. Here $\overline M$ is called the base, $\widetilde M$ is called the fiber and $f$ is called warping function of the warped product $M = \overline M \times_f \widetilde M$. If $f=1$, then the warped product reduces to semi-Riemannian product. We denote $\Gamma^i_{jk}$, $R_{ijkl}$, $S_{ij}$ and $\kappa$ as the components of Levi-Civita connection $\nabla$, the Riemann-Christoffel curvature tensor $R$, Ricci tensor $S$ and the scalar curvature of $(M, g)$ respectively. Moreover we consider that, when $\Omega$ is a quantity formed with respect to $g$, we denote by $\overline \Omega$ and $\widetilde \Omega$, the similar quantities formed with respect to $\overline g$ and $\widetilde g$ respectively.\\
%=====================================================
\indent Then the non-zero local components of Levi-Civita connection $\nabla$ of $(M,g)$ are given by
\be\label{eq2.2}
\Gamma^a_{bc}=\overline{\Gamma}^a_{bc},\,\,\,\, \Gamma^\alpha_{\beta \gamma}=\widetilde{\Gamma}^\alpha_{\beta \gamma},\,\,\,\,\,\,\,\Gamma^a_{\beta \gamma}=-\frac{1}{2}\overline{g}^{ab}f_{b} \widetilde{g}_{\beta \gamma},\,\,\ \ \Gamma^\alpha_{ a \beta }=\frac{1}{2f}f_{a}\delta^{\alpha}_{\beta},
\ee
where $f_{a}=\partial_{a} f=\frac{\partial f}{\partial x^{a}}$.\\
%----------------------------------------------------
The local components $R_{hijk}=g_{hl}R^{l}_{ijk}=g_{hl}(\partial_{k}\Gamma^{l}_{ij}-\partial_{j}\Gamma^{l}_{ik}+  \Gamma^{m}_{ij}\Gamma^{l}_{mk}-\Gamma^{m}_{ik}\Gamma^{l}_{mj}),\,\,\, \partial_{k}=\frac{\partial}{\partial x^{k}},$
of the Riemann-Christoffel curvature tensor $R$ of $(M,g)$ which may not vanish identically are the following:
\be\label{R}
R_{abcd} = \overline{R}_{abcd},\,\,\,\, R_{a\alpha b\beta}=f T_{ab}\widetilde{g}_{\alpha \beta},\,\,\,R_{\alpha \beta \gamma \delta} = f\widetilde{R}_{\alpha \beta \gamma \delta} - f^2 P \widetilde{G}_{\alpha \beta \gamma \delta}, 
\ee
where $G_{ijkl} = g_{il}g_{jk}-g_{ik}g_{jl}$ are the components of Gaussian curvature and
$$T_{ab} = -\frac{1}{2f}(\nabla_b f_a - \frac{1}{2f}f_a f_b), \ \ \ \ tr(T) = g^{ab}T_{ab},$$
$$Q = f((n-p-1)P -tr(T)), \ \ \ \ P = \frac{1}{4f^2}g^{ab}f_a f_b.$$
%----------------------------------------------------
Again, the non-zero local components of the Ricci tensor $S_{jk} = g^{il}R_{ijkl}$ of $(M, g)$ are given by
\be\label{eq2.4}
S_{ab}=\overline{S}_{ab}-(n-p)T_{ab},\,\,\,\, S_{\alpha \beta}=\widetilde{ S}_{\alpha \beta} + Q \widetilde{g}_{\alpha \beta}.  
\ee
The scalar curvature $\kappa$ of $(M, g)$ is given by
\be\label{eq2.5}
\kappa=\overline{\kappa}+\frac{\widetilde{\kappa}}{f}-(n-p)[(n-p-1)P - 2 \; tr(T)].  
\ee
For more detail about warped product components of basic tensors we refer the readers to see \cite{Hotl04}, \cite{SK12} and also references therein.\\
%=============================================================================================
\indent Now from the above we can easily calculate the components of various necessary tensors of a warped product manifold in terms of its base and fiber components. The non-zero components of Ricci tensor of level 2 are 
\bea
\left\{
\begin{array}{l}
(i) S^2_{ab} = \overline{S}^2_{ab}+ (n-p)(\overline S\cdot T)_{ab}+(n-p)^2 T^2_{ab},\\
(ii) S^2_{\alpha\beta} = \frac{1}{f}[\widetilde S^2_{\alpha\beta} + 2Q \widetilde S_{\alpha\beta} + Q^2 \widetilde g_{\alpha\beta}].
\end{array}
\right.
\eea
%----------------------------------------------------
The non-zero components of $(g\wedge g)$ are 
\bea\label{gg}
\left\{
\begin{array}{l}
(i) (g\wedge g)_{abcd} = (\overline g \wedge \overline g)_{abcd},\\
(ii) (g\wedge g)_{a\alpha b \beta} = -2 f \overline g_{ab}\widetilde g_{\alpha\beta},\\
(iii) (g\wedge g)_{\alpha\beta\gamma\delta} = f^2 (\widetilde g \wedge \widetilde g)_{\alpha\beta\gamma\delta}.
\end{array}
\right.
\eea
%----------------------------------------------------
The non-zero components of $(g\wedge S)$ are 
\bea\label{gs}
\left\{
\begin{array}{l}
(i) (g\wedge S)_{abcd} = (\overline g \wedge \overline S)_{abcd} - (n-p) (\overline g \wedge T)_{abcd},\\
(ii) (g\wedge S)_{a\alpha b \beta} = -\overline g_{ab}(\widetilde S_{\alpha\beta}+Q\widetilde g_{\alpha\beta})-f \widetilde g_{\alpha\beta}(\overline S_{ab}-(n-p)T_{ab}),\\
(iii) (g\wedge S)_{\alpha\beta\gamma\delta} = f (\widetilde g \wedge \widetilde S)_{\alpha\beta\gamma\delta} 
+ 2 f Q \widetilde G_{\alpha\beta\gamma\delta}.
\end{array}
\right.
\eea
%----------------------------------------------------
The non-zero components of $(S\wedge S)$ are 
\bea\label{ss}
\left\{
\begin{array}{l}
(i) (S\wedge S)_{abcd} = (\overline S \wedge \overline S)_{abcd} - (n-p) (\overline S \wedge T)_{abcd} + (n-p)^2 (T \wedge T)_{abcd},\\
(ii) (S\wedge S)_{a\alpha b \beta} = -2(\widetilde S_{\alpha\beta}+Q\widetilde g_{\alpha\beta})(\overline S_{ab}-(n-p)T_{ab}),\\
(iii) (S\wedge S)_{\alpha\beta\gamma\delta} = (\widetilde S \wedge \widetilde S)_{\alpha\beta\gamma\delta} 
+ Q (\widetilde S \wedge \widetilde g)_{\alpha\beta\gamma\delta} + Q^2 (\widetilde g \wedge \widetilde g)_{\alpha\beta\gamma\delta}.
\end{array}
\right.
\eea
%----------------------------------------------------
The non-zero components of $(g\wedge S^2)$ are 
\bea\label{gsq}
\left\{
\begin{array}{l}
(i) (g\wedge S^2)_{abcd} = (\overline g \wedge \overline S^2)_{abcd} 
+ (n-p) (\overline g \wedge (\overline S\cdot T))_{abcd} + (n-p)^2 (\overline g \wedge T^2)_{abcd},\\
(ii) (g\wedge S^2)_{a\alpha b \beta} = 
-\frac{1}{f}\overline g_{ab}(\widetilde S^2_{\alpha\beta}+ 2 Q\widetilde S_{\alpha\beta} + Q^2 \widetilde g_{\alpha\beta})\\
\hspace{1.3in} -f\widetilde g_{\alpha\beta}(\overline S^2_{ab}+(n-p)\overline S\cdot T_{ab} + (n-p)^2 T^2_{ab}),\\
(iii) (g\wedge S^2)_{\alpha\beta\gamma\delta} = 
(\widetilde g \wedge \widetilde S^2)_{\alpha\beta\gamma\delta} 
+ 2 Q (\widetilde g \wedge \widetilde S)_{\alpha\beta\gamma\delta} 
+ Q^2 (\widetilde g \wedge \widetilde g)_{\alpha\beta\gamma\delta}.
\end{array}
\right.
\eea
%----------------------------------------------------
The non-zero components of $(S\wedge S^2)$ are 
\bea\label{ssq}
\left\{
\begin{array}{l}
(i) (S\wedge S^2)_{abcd} = (\overline S \wedge \overline S^2)_{abcd} 
+ (n-p) (\overline S \wedge (\overline S\cdot T))_{abcd}\\
\hspace{1.2in}+ (n-p)^2 (\overline S \wedge T^2)_{abcd}-(n-p)(\overline S^2 \wedge T)_{abcd}\\
\hspace{1.2in} -(n-p)^2(T \wedge (\overline S\cdot T))_{abcd} + (n-p)^3 (T \wedge T^2)_{abcd},\\
(ii) (S\wedge S^2)_{a\alpha b \beta} = 
-\frac{1}{f}(\overline S_{ab}-(n-p)T_{ab})(\widetilde S^2_{\alpha\beta}+ 2 Q\widetilde S_{\alpha\beta} + Q^2 \widetilde g_{\alpha\beta})\\
\hspace{1.3in} -(\overline S^2_{ab}+(n-p)(\overline S\cdot T)_{ab} + (n-p)^2 T^2_{ab})(\widetilde S_{\alpha\beta} + Q\widetilde g_{\alpha\beta}),\\
(iii) (S\wedge S^2)_{\alpha\beta\gamma\delta} = 
\frac{1}{f}[(\widetilde S \wedge \widetilde S^2)_{\alpha\beta\gamma\delta} 
+ 4 Q (\widetilde S \wedge \widetilde S)_{\alpha\beta\gamma\delta} 
+ Q^2 (\widetilde S \wedge \widetilde g)_{\alpha\beta\gamma\delta}\\
\hspace{1.2in}+ Q (\widetilde g \wedge \widetilde S^2)_{\alpha\beta\gamma\delta}
+ 2 Q^2 (\widetilde g \wedge \widetilde S)_{\alpha\beta\gamma\delta}
+ 2Q^3 (\widetilde g \wedge \widetilde g)_{\alpha\beta\gamma\delta}].
\end{array}
\right.
\eea
%----------------------------------------------------
The non-zero components of $(S^2\wedge S^2)$ are 
\bea\label{sqsq}
\left\{
\begin{array}{l}
(i) (S^2\wedge S^2)_{abcd} = (\overline S^2 \wedge \overline S^2)_{abcd} + (n-p)^2 ((\overline S\cdot T) \wedge (\overline S\cdot T))_{abcd}\\
\hspace{1.2in}+(n-p)^2 (T^2 \wedge T^2)_{abcd}+ 2(n-p) (\overline S^2 \wedge (\overline S\cdot T))_{abcd}\\
\hspace{1.2in}+2(n-p)^3((\overline S\cdot T^2) \wedge T^2)_{abcd}+ 2(n-p)^3 (\overline S^2 \wedge T^2)_{abcd},\\
(ii) (S^2\wedge S^2)_{a\alpha b \beta} = 
-\frac{2}{f}(\overline S^2_{ab}+(n-p)(\overline S\cdot T)_{ab} + (n-p)^2 T^2_{ab})\\
\hspace{1.7in}(\widetilde S^2_{\alpha\beta}+ 2 Q\widetilde S_{\alpha\beta} + Q^2 \widetilde g_{\alpha\beta}),\\
(iii) (S\wedge S^2)_{\alpha\beta\gamma\delta} = 
\frac{1}{f^2}[(\widetilde S^2 \wedge \widetilde S^2)_{\alpha\beta\gamma\delta} 
+ 4 Q^2 (\widetilde S \wedge \widetilde S)_{\alpha\beta\gamma\delta} 
+ Q^4 (\widetilde g \wedge \widetilde g)_{\alpha\beta\gamma\delta}\\
\hspace{1.2in}+ 4 Q (\widetilde S^2 \wedge \widetilde S)_{\alpha\beta\gamma\delta}
+ 2 Q^2 (\widetilde g \wedge \widetilde S^2)_{\alpha\beta\gamma\delta}
+ 4 Q^3 (\widetilde g \wedge \widetilde S)_{\alpha\beta\gamma\delta}].
\end{array}
\right.
\eea
%===============================================================================================================
\indent From above we see that the components of $g\wedge g$, $g\wedge S$, $S\wedge S$, $g\wedge S^2$, $S\wedge S^2$ and $S^2\wedge S^2$ are in a quadratic form of wedge product for base and fiber part and quadratic form of the product for the mixed part. So each of them can be expressed by a matrix. For example, $(g\wedge S)_{abcd}$,  $(g\wedge S)_{a\alpha c\beta}$ and $(g\wedge S)_{\alpha\beta\gamma\delta}$ can respectively be expressed as:\\
%===================================================
$\left(
\begin{array}{c}
 \overline g \\
 \overline S \\
 \overline S^2 \\
 T \\
 T^2 \\
 \overline S\cdot T
\end{array}
\right)^t$
%---------------
$\left(
\begin{array}{cccccc}
 0 & \frac{1}{2} & 0 & \frac{p-n}{2} & 0 & 0 \\
 \frac{1}{2} & 0 & 0 & 0 & 0 & 0 \\
 0 & 0 & 0 & 0 & 0 & 0 \\
 \frac{p-n}{2} & 0 & 0 & 0 & 0 & 0 \\
 0 & 0 & 0 & 0 & 0 & 0 \\
 0 & 0 & 0 & 0 & 0 & 0
\end{array}
\right)$ 
%------------------
$\wedge$
%------------------
$\left(
\begin{array}{c}
 \overline g \\
 \overline S \\
 \overline S^2 \\
 T \\
 T^2 \\
 \overline S\cdot T
\end{array}
\right)_{abcd}$ or 
%=-=-=-=-=-=-=-=-=-=-=-=-=-=-=
$\begin{array}{|c|c|c|c|c|c|c|}\hline
 \land  & \overline g & \overline S & \overline S^2 & T & T^2 & \overline S\cdot T \\\hline
 \overline g & 0 & \frac{1}{2} & 0 & \frac{p-n}{2} & 0 & 0 \\\hline
 \overline S & \frac{1}{2} & 0 & 0 & 0 & 0 & 0 \\\hline
 \overline S^2 & 0 & 0 & 0 & 0 & 0 & 0 \\\hline
 T & \frac{p-n}{2} & 0 & 0 & 0 & 0 & 0 \\\hline
 T^2 & 0 & 0 & 0 & 0 & 0 & 0 \\\hline
 \overline S\cdot T & 0 & 0 & 0 & 0 & 0 & 0\\\hline
\end{array}$\\
%==================================================
%
$\left(
\begin{array}{c}
 \overline g \\
 \overline S \\
 \overline S^2 \\
 T \\
 T^2 \\
 \overline S\cdot T
\end{array}
\right)^t_{ab}$
%----------------------
$\left(
\begin{array}{ccc}
 -Q & -1 & 0 \\
 -f & 0 & 0 \\
 0 & 0 & 0 \\
 f (p-n) & 0 & 0 \\
 0 & 0 & 0 \\
 0 & 0 & 0
\end{array}
\right)$
%-----------------
$\left(
\begin{array}{c}
 \widetilde g \\
 \widetilde S \\
 \widetilde S^2
\end{array}
\right)_{\alpha\beta}$ or
%=-=-=-=-=-=-=-=-=-=-==-=-=-=-=-=-=-=
$\begin{array}{|c|c|c|c|}\hline
   & \widetilde{g} & \widetilde{S} & \widetilde S^2 \\\hline
 \overline{g} & -Q & -1 & 0 \\\hline
 \overline{S} & -f & 0 & 0 \\\hline
 \overline S^2 & 0 & 0 & 0 \\\hline
 T & f (p-n) & 0 & 0 \\\hline
 T^2 & 0 & 0 & 0 \\\hline
 \overline S\cdot T & 0 & 0 & 0\\\hline
\end{array}$\\
%====================================
%
$\left(
\begin{array}{c}
 \widetilde g \\
 \widetilde S \\
 \widetilde S^2
\end{array}
\right)^t$
%----------------
$\left(
\begin{array}{ccc}
 f Q & \frac{f}{2} & 0 \\
 \frac{f}{2} & 0 & 0 \\
 0 & 0 & 0
\end{array}
\right)$
%--------------------
$\wedge$
%--------------------
$\left(
\begin{array}{c}
 \widetilde g \\
 \widetilde S \\
 \widetilde S^2
\end{array}
\right)_{\alpha\beta\gamma\delta}$ or
%=-=-=-=-=-=-==-=-=-=-=-=-=-=-=-=-=-=-=-=
$\begin{array}{|c|c|c|c|}\hline
 \land  & \widetilde g & \widetilde S & \widetilde S^2 \\\hline
 \widetilde g & f Q & \frac{f}{2} & 0 \\\hline
 \widetilde S & \frac{f}{2} & 0 & 0 \\\hline
 \widetilde S^2 & 0 & 0 & 0\\\hline
\end{array}$\\
%=========================================
Similarly we can get the matrix representation for the other tensors of $g\wedge g$, $S\wedge S$, $g\wedge S^2$, $S\wedge S^2$ and $S^2\wedge S^2$.
%%%%%%%%%%%%%%%%%%%%%%%%%%%%%%%%%%%%%%%%%%%%%%%%%%%%%%%%%%%%%%%%%%%%%%%%%%%%%%%%%%%%%%%%%%%%%%%%%%%%%%%%%%%%%%
%                                5. Warped Product $GRT$
%%%%%%%%%%%%%%%%%%%%%%%%%%%%%%%%%%%%%%%%%%%%%%%%%%%%%%%%%%%%%%%%%%%%%%%%%%%%%%%%%%%%%%%%%%%%%%%%%%%%%%%%%%%%%%%%
\section{\bf Warped Product generalized Roter-type manifolds}\label{main}
%****************************************************************************
\begin{thm}\label{th5.1}
Let $M^n = \overline M^p \times_f \widetilde M^{n-p}$ be a warped product manifold. Then $M$ is a generalized Roter-type with
\be\label{eq5.1}
R= L_1 g\wedge g + L_2 g\wedge S + L_3 S\wedge S + L_4 g\wedge S^2 + L_5 S\wedge S^2 + L_6 S^2\wedge S^2
\ee
if and only if\\
%-------------------------------
$(i)$ the Riemann-Christoffel curvature tensor $\overline R$ of $\overline M$ can be expressed as\\
$\begin{array}{|c|c|c|c|c|c|c|}\hline
 \land  & \overline g & \overline S & \overline S^2 & T & T^2 & \overline S\cdot T \\\hline
 \overline g & L_1 & \frac{L_2}{2} & \frac{L_4}{2} & \frac{1}{2} L_2 (p-n) & \frac{1}{2} L_4 (n-p)^2 & \frac{1}{2} L_4 (n-p) \\\hline
 \overline S & \frac{L_2}{2} & L_4 & \frac{L_5}{2} & L_4 (p-n) & \frac{1}{2} L_5 (n-p)^2 & \frac{1}{2} L_5 (n-p) \\\hline
 \overline S^2 & \frac{L_4}{2} & \frac{L_5}{2} & L_6 & \frac{1}{2} L_5 (p-n) & L_6 (n-p)^2 & L_6 (n-p) \\\hline
 T & \frac{1}{2} L_2 (p-n) & L_3 (p-n) & \frac{1}{2} L_5 (p-n) & L_3 (n-p)^2 & -\frac{1}{2} L_5 (n-p)^3 & -\frac{1}{2} L_5(n-p)^2 \\\hline
 T^2 & \frac{1}{2} L_4 (n-p)^2 & \frac{1}{2} L_5 (n-p)^2 & L_6 (n-p)^2 & -\frac{1}{2} L_5 (n-p)^3 & L_6 (n-p)^4 & L_6(n-p)^3 \\\hline
 \overline S\cdot T & \frac{1}{2} L_4 (n-p) & \frac{1}{2} L_5 (n-p) & L_6 (n-p) & -\frac{1}{2} L_5 (n-p)^2 & L_6 (n-p)^3 & L_6 (n-p)^2\\\hline
\end{array}$\\
%-------------------------------
$(ii)$ the following expression is equal to zero\\
$\begin{array}{|c|c|c|c|}\hline
   & \widetilde{g} & \widetilde{S} & \widetilde S^2 \\\hline
 \overline{g} & -\frac{L_4 Q^2}{f}-L_2 Q-2 f L_1 & -L_2-\frac{2 L_4 Q}{f} & -\frac{L_4}{f} \\\hline
 \overline{S} & -\frac{L_5 Q^2}{f}-2 L_3 Q-f L_2 & -\frac{2 (f L_3+L_5 Q)}{f} & -\frac{L_5}{f} \\\hline
 \overline S^2 & -\frac{2 L_6 Q^2}{f}-L_5 Q-f L_4 & -L_5-\frac{4 L_6 Q}{f} & -\frac{2 L_6}{f} \\\hline
 T & f (L_2 (p-n)-1)+\frac{L_5 Q^2 (p-n)}{f}+2 L_3 Q (p-n) & -\frac{2 (n-p) (f L_3+L_5 Q)}{f} & \frac{L_5 (p-n)}{f} \\\hline
 T^2 & -\frac{(n-p)^2 \left(L_4 f^2+L_5 Q f+2 L_6 Q^2\right)}{f} & -\frac{(n-p)^2 (f L_5+4 L_6 Q)}{f} & -\frac{2 L_6 (n-p)^2}{f} \\\hline
 \overline S\cdot T & -\frac{(n-p) \left(L_4 f^2+L_5 Q f+2 L_6 Q^2\right)}{f} & -\frac{(n-p) (f L_5+4 L_6 Q)}{f} & \frac{2 L_6 (p-n)}{f}\\\hline
\end{array}$\\
$(iii)$ $f \widetilde R$, $\widetilde R$ be the Riemann-Christoffel curvature tensor of $\widetilde M$, can be expressed as\\
$\begin{array}{|c|c|c|c|}\hline
 \land  & \widetilde g & \widetilde S & \widetilde S^2 \\\hline

 \widetilde g & \frac{L_6 Q^4}{f^2}+\frac{L_5 Q^3}{f}+(L_3+L_4) Q^2+ &
 (L_3+L_4) Q+& \frac{1}{2} \left(L_4+\frac{Q (f L_5+2 L_6 Q)}{f^2}\right) \\

 $ $ & f L_2 Q+f^2 L_1 - \frac{f P}{2} & \frac{L_2 f^3+3 L_5 Q^2 f+ 4 L_6 Q^3}{2 f^2} & $ $ \\\hline

 \widetilde S & (L_3+L_4) Q+\frac{L_2 f^3+3 L_5 Q^2 f+4 L_6 Q^3}{2 f^2} &
 L_3+\frac{2 Q (f L_5+2 L_6 Q)}{f^2} &
 \frac{f L_5+4 L_6 Q}{2 f^2} \\\hline

 \widetilde S^2 & \frac{1}{2} \left(L_4+\frac{Q (f L_5+2 L_6 Q)}{f^2}\right) &
 \frac{f L_5+4 L_6 Q}{2 f^2} &
 \frac{L_6}{f^2}\\\hline
\end{array}$
\end{thm}
\noindent {\bf Proof:} By putting the values of $R$, $g\wedge g$, $g\wedge S$, $S\wedge S$, $g\wedge S^2$, $S\wedge S^2$, $S^2\wedge S^2$ for all possible cases of their components from (\ref{R}) and (\ref{gg}) to (\ref{sqsq}) in the generalized Roter type condition (\ref{eq5.1}) we get our assertion easily.\\
%%%%%%%%%%%%%%%%%%%%%%%%%%%%%%%%%%%%%%%%%%%%%%%%%%
\indent From the above we can conclude on the base and fiber part of the warped product generalized Roter type manifold as follows:
\begin{cor}\label{cor5.1}
Let $M^n = \overline M^p \times_f \widetilde M^{n-p}$ be a warped product manifold satisfying generalized Roter-type condition as
$$R= L_1 S\wedge S + L_2 S\wedge S^2 + L_3 g\wedge S + L_4 g\wedge S^2 + L_5 g\wedge g + L_6 S^2\wedge S^2.$$
Then, 
%------------------------
$($i$)$ the fiber $\widetilde M$ is generalized Roter type.\\
$($ii$)$ the fiber $\widetilde M$ is Roter type if $J_1 \neq 0$,\\
where $J_1 = -\frac{L_5 \overline{r}+2 L_6 \left((n-p) \left(tr(T^2)(n-p)+tr(\overline{S}\cdot T)\right)+\overline{r^{(2)}}\right)+L_4 p+L_5 tr(T) (n-p)}{f}$. Moreover in this case fiber satisfies $Ein(2)$ condition.\\
$($iii$)$ the fiber $\widetilde M$ is conformally flat if $J_1 \neq 0$ and
$$\frac{(J_2)^2 L_6}{f^2 (J_1)^2}+\frac{J_2 (f L_5+4 L_6 Q)}{f^2 J_1}+\frac{2 Q (f L_5+2 L_6 Q)}{f^2}+L_3 = 0,$$
where 
\beb
J_2 = &-& \frac{1}{f}\Big[2 \overline{r} (f L_3+L_5 Q)+(f L_5+4 L_6 Q) \left((n-p) \left(tr(T^2)(n-p)+tr(\overline{S}\cdot T)\right)+\overline{r^{(2)}}\right)\\
&+&p (f L_2+2 L_4 Q)+2 tr(T) (n-p) (f L_3+L_5 Q)\Big].
\eeb
$($iv$)$ the fiber $\widetilde M$ is of constant curvature if $J_1 = 0$ and $J_2\neq 0$. Moreover in this case the fiber is Einstein.\\
$($v$)$ the base $\overline M$ is generalized Roter type if $T$, $\overline g$ and $\overline S$ are linearly dependent with non-zero coefficient of $T$.
\end{cor}
%%%%%%%%%%%%%%%%%%%%%%%%%%%%%%%%%%%%%%%%%%%%%%%%%%
From Theorem \ref{th5.1} we can easily get the necessary and sufficient condition for a warped product manifold to be Roter type.
%==============================
\begin{cor}\label{cor5.2}
Let $M= \overline M\times_f \widetilde M$ be a non-flat warped product manifold. Then $M$ is a Roter type manifold with
$$R = N_1 g\wedge g + N_2 g\wedge S + N_3 S\wedge S$$
if and only if\\
%--------------------
(i) the Rieman-Christoffel curvature tensor $\overline R$ of $\overline M$ can be expressed as\\
$\begin{array}{|c|c|c|c|}\hline
 \land  & \overline g & \overline S & T \\\hline
 \overline g & N_1 & \frac{N_2}{2} & \frac{1}{2} N_2 (p-n) \\\hline
 \overline S & \frac{N_2}{2} & N_3 & N_3 (p-n) \\\hline
 T & \frac{1}{2} N_2 (p-n) & N_3 (p-n) & N_3 (n-p)^2\\\hline
\end{array}$\\
%------------------
(ii) the following expression is equal to zero\\
$\begin{array}{|c|c|c|c|}\hline
   & \overline{g} & \overline{S} & T \\\hline
 \widetilde{g} & -2 f N_1-N_2 Q & -f N_2-2 N_3 Q & -f (1+N_2 (n-p))-2 N_3 Q (n-p) \\\hline
 \widetilde{S} & -N_2 & -2 N_3 & 2 N_3 (p-n)\\\hline
\end{array}$\\
%-----------------
(iii) $f \widetilde R$, $\widetilde R$ be the Riemann-Christoffel curvature tensor of $\widetilde M$, can be expressed as\\
$\begin{array}{|c|c|c|}\hline
 \land  & \widetilde g & \widetilde S \\\hline
 \widetilde g & N_1 f^2+N_2 Q f+N_3 Q^2 -\frac{f P}{2}& \frac{f N_2}{2}+N_3 Q \\\hline
 \widetilde S & \frac{f N_2}{2}+N_3 Q & N_3\\\hline
\end{array}$
\end{cor}
%============================
As similar to the Corollary \ref{cor5.1} we can conclude the following on the base and fiber part of a Roter type warped product manifold.
\begin{cor}\label{cor5.3}
Let $M= \overline M\times_f \widetilde M$ be a non-flat warped product manifold. Then $M$ is a Roter type manifold with
$$R = N_1 g\wedge g + N_2 g\wedge S + N_3 S\wedge S.$$
Then
$($i$)$ fiber is of Roter type.\\
$($ii$)$ fiber is conformally flat if $M$ is conformally flat.\\
$($iii$)$ fiber is of constant curvature if $-2(n-p) N_3 tr(T)-N_2 p-2 N_3 r \neq 0$, and in this case fiber is an Einstein manifold.\\
$($iv$)$ base is of Roter type if $T$, $\overline g$ and $\overline S$ are linearly dependent with non-zero coefficient of $T$.
\end{cor}
%%%%%%%%%%%%%%%%%%%%%%%%%%%%%%%%%%%%%%%%%%%%%%%%%%
Now we can easily deduce the necessary and sufficient condition for a warped product manifold to be conformally flat, as follows:
%==============================
\begin{cor}\label{cor5.4}
Let $M= \overline M\times_f \widetilde M$ be a non-flat warped product manifold. Then $M$ is conformally flat if and only if\\
%--------------------
(i) $\overline R = \frac{\kappa}{(n-2) (n-1)}\overline g\wedge \overline g +\frac{1}{n-2}\overline g\wedge \overline S -\frac{n-p}{n-2}\overline g\wedge T$\\
%------------------
(ii) $\left[-\frac{2 f \kappa}{(n-2) (n-1)}-\frac{Q}{n-2}\right]\overline g_{ab}\widetilde g_{\alpha\beta}
-\frac{1}{n-2}\overline g_{ab}\widetilde S_{\alpha\beta}
-\frac{f}{n-2}\overline S_{ab}\widetilde g_{\alpha\beta}
-f \left(\frac{n-p}{n-2}+1\right)T_{ab}\widetilde g_{\alpha\beta} = 0$\\
%-----------------
(iii) $\widetilde R = \left[\frac{f \kappa}{(n-2) (n-1)}+\frac{Q}{n-2}-\frac{1}{2} P\right] \widetilde g\wedge \widetilde g +\frac{1}{(n-2)}\widetilde g\wedge \widetilde S$.
\end{cor}
%----------------
\noindent \textbf{Proof:} The result follows from Corollary \ref{cor5.2} by using $N_1 = \frac{r}{(n-1)(n-2)}$, $N_2 = \frac{1}{n-2}$ and $N_3=0$.\\
%============================
From above we can state the following:
\begin{cor}
In a conformally flat warped product manifold\\
(i) fiber is conformally flat.\\
(ii) fiber is of constant curvature if and only if it is Einstein.\\
(iii) fiber is of quasi-constant curvature if and only if it is quasi-Einstein.\\
$($iv$)$ base is of conformally flat if $T$, $\overline g$ and $\overline S$ are linearly dependent with non-zero coefficient of $T$.
\end{cor}
%%%%%%%%%%%%%%%%%%%%%%%%%%%%%%%%======================================%%%%%%%%%%%%%%%%%%%%%%%%%%%
We now discuss about the decomposable or product semi-Riemannian manifold satisfying some generalized Roter type conditions. We know that semi-Riemannian product are some special case of warped product manifold, where the warping function $f$ is identically 1. Then we have
$$T =0, \ \ P =0 \ \mbox{and } \ Q=0.$$
Thus applying these values in (\ref{R}) to (\ref{sqsq}) we get the non-zero components of $R$, $S$, $r$, $S^2$, $g\wedge g$, $g\wedge S$, $S\wedge S$, $g\wedge S^2$, $S\wedge S^2$ and $S^2\wedge S^2$. Now from Theorem \ref{th5.1} we can state the following:
%======================
\begin{cor}\label{cor5.5}
Let $M^n = \overline M^p \times \widetilde M^{n-p}$ be a product manifold. Then $M$ is a generalized Roter-type with
\be
R= L_1 g\wedge g + L_2 g\wedge S + L_3 S\wedge S + L_4 g\wedge S^2 + L_5 S\wedge S^2 + L_6 S^2\wedge S^2
\ee
if and only if\\
%-------------------------------
$(i)$ the Riemann-Christoffel curvature tensor $\overline R$ of $\overline M$ and $\widetilde R$ of $\widetilde M$ can respectively be expressed as\\
$\begin{array}{|c|c|c|c|}\hline
 \land  & \overline g & \overline S & \overline S^2 \\\hline
 \overline g & L_1 & \frac{L_2}{2} & \frac{L_4}{2} \\\hline
 \overline S & \frac{L_2}{2} & L_4 & \frac{L_5}{2} \\\hline
 \overline S^2 & \frac{L_4}{2} & \frac{L_5}{2} & L_6 \\\hline
\end{array}$ \ \ \ and \ \ \ 
%-------------------------------
$\begin{array}{|c|c|c|c|}\hline
 \land  & \widetilde g & \widetilde S & \widetilde S^2 \\\hline

 \widetilde g & L_1 & \frac{L_2}{2} & \frac{L_4}{2} \\\hline

 \widetilde S & \frac{L_2}{2} & L_3 & \frac{L_5}{2} \\\hline

 \widetilde S^2 & \frac{L_4}{2} & \frac{L_5}{2} & L_6\\\hline
\end{array}$\\
%-------------------------------
$(ii)$ the following expression is equal to zero\\
$\begin{array}{|c|c|c|c|}\hline
   & \widetilde{g} & \widetilde{S} & \widetilde S^2 \\\hline
 \overline{g} & -2 L_1 & -L_2 & -L_4 \\\hline
 \overline{S} & -L_2 & -2 L_3 & -L_5 \\\hline
 \overline S^2 & -L_4 & -L_5 & -2 L_6 \\\hline
\end{array}$
\end{cor}
%====================================================
\noindent\textbf{Note:} From the above corollary we can get a necessary and sufficient condition for a product manifold to be Roter type by taking $L_4 = L_5 = L_6 = 0$, and conformally flat by taking  $L_3 = L_4 = L_5 = L_6 = 0$, $L_1 = \frac{r}{2(n-1)(n-2)}$ and $L_2 = \frac{1}{n-2}$. Again from the above results we see that the decompositions of a semi-Riemannian product generalized Roter type manifold are also generalized Roter type manifold but the converse is not necessarily true, in general (see Example 5.1). We also note that the same case arises for Roter type and conformally flat manifolds also (see Example 5.1).
%%%%%%%%%%%%%%%%%%%%%%%%%%%%%%%%%%%%%%%%%%%%%%%%%%%%%%%%%%%%%%%%%%%%%%%%%%%%%%%%%%%%%%
%                                              Examples
%%%%%%%%%%%%%%%%%%%%%%%%%%%%%%%%%%%%%%%%%%%%%%%%%%%%%%%%%%%%%%%%%%%%%%%%%%%%%%%%%%%%%%
\section{\bf Examples}\label{exam}
%=========================
%****************************************************************************
\textbf{Example 5.1:} Consider the warped product $M = \overline M\times_f \widetilde M$, where $\overline M$ is an open interval of $\mathbb R$ with usual metric $\overline g = (dx^1)^2$ in local coordinate $x^1$ and $\widetilde M$ is a 4-dimensional manifold equipped with a semi-Riemannian metric
$$\widetilde g = (dx^2)^2+h(dx^3)^2+h(dx^4)^2+h\psi(dx^5)^2$$
in local coordinates $(x^2, x^3, x^4, x^5)$, where the warping function $f$ is a function of $x^1$ and the functions $h$ and $\psi$ are non-zero functions of $x^2$ and $x^3$ respectively. We can easily evaluate the local components of necessary tensors of $\widetilde M$. The local non-zero components of the Riemann-Christoffel curvature tensor $\widetilde R$ and the Ricci tensor $\widetilde S$ of $\widetilde M$ upto symmetry are
$$\psi \widetilde R_{1212}= \psi \widetilde R_{1313}= \widetilde R_{1414}=\psi  \frac{\left(\left(h'\right)^2-2 h h''\right)}{4 h},$$
$$\psi \widetilde R_{2323}= \widetilde R_{3434}=-\frac{\psi}{4}  \left(h'\right)^2,$$
$$\widetilde R_{2424}=\frac{1}{4} \left(-\psi  \left(h'\right)^2-2 h \psi ''+\frac{h \left(\psi '\right)^2}{\psi }\right)$$
and
$$\widetilde S_{11}=\frac{3 \left(2 h h''-\left(h'\right)^2\right)}{4 h^2},$$
$$\widetilde S_{22}=\frac{1}{4} \left(2 h''+\frac{\left(h'\right)^2}{h}-\frac{\left(\psi '\right)^2-2 \psi  \psi ''}{\psi ^2}\right),$$
$$\widetilde S_{33}=\frac{2 h h''+\left(h'\right)^2}{4 h},$$
$$\widetilde S_{44}=\frac{1}{4} \left(2 \left(\psi  h''+\psi ''\right)+\frac{\psi  \left(h'\right)^2}{h}-\frac{\left(\psi '\right)^2}{\psi }\right).$$
Then we can easily check that this manifold is of generalized Roter type and satisfies the $Ein(3)$ condition. Again if\\
(i) $(h')^2-h h'' =0$, i.e., $h=c_1 e^{c_2 x^2}$, then it satisfies the $Ein(2)$ condition and thus becomes Roter type;\\
(ii) $(\psi')^2 - 2 \psi \psi'' =0$, i.e., $\psi =\frac{(c_1 x^3 +2 c_2)^2}{4 c_2}$, then it becomes a manifold of constant curvature,\\
where $c_1$ and $c_2$ are arbitrary constants.\\
%------------------
\indent Thus by a straightforward calculation we can evaluate the components of various necessary tensors corresponding to $M$. The non-zero local components of the Riemann-Christoffel curvature tensor $R$ and the Ricci tensor $S$ of $M$ upto symmetry are
%-------------------------
$$h \psi R_{1212}= \psi R_{1313}= \psi R_{1414}= R_{1515}=h \psi \frac{\left(f'\right)^2-2 f f''}{4 f},$$
$$\psi R_{2323}= \psi R_{2424} = R_{2525} = \frac{\psi}{4} \left(-h \left(f'\right)^2-2 f h''+\frac{f \left(h'\right)^2}{h}\right),$$
$$\psi R_{3434}= R_{4545}=-\frac{\psi}{4} \left(h^2 \left(f'\right)^2+f \left(h'\right)^2\right)$$
$$R_{3535}=\frac{1}{4} \left[f \left(-\psi  \left(h'\right)^2-2 h \psi ''+\frac{h \left(\psi '\right)^2}{\psi }\right)-h^2 \psi  \left(f'\right)^2\right].$$
and
$$S_{11}=-\frac{\left(f'\right)^2-2 f f''}{f^2},$$
$$S_{22}=\frac{1}{4} \left(2 f''+\frac{2 \left(f'\right)^2}{f}+\frac{6 h h''-3 \left(h'\right)^2}{h^2}\right),$$
$$\psi S_{33} = S_{55} = \frac{1}{4} \left(2 h \psi  f''+\frac{2 h \psi  \left(f'\right)^2}{f}+2 \psi  h''+\frac{\psi  \left(h'\right)^2}{h}+2 \psi ''-\frac{\left(\psi '\right)^2}{\psi }\right),$$
$$S_{44}=\frac{1}{4} \left(2 \left(h f''+h''\right)+\frac{2 h \left(f'\right)^2}{f}+\frac{\left(h'\right)^2}{h}\right).$$
From these we can easily calculate the local components of $S^2$, $S^3$, $S^4$ and also the local components of $G$, $g\wedge S$, $S\wedge S$, $g\wedge S^2$, $S\wedge S^2$ and $S^2\wedge S^2$. We observe that for any $f$, $h$ and $\psi$, the manifold is $Ein(4)$ but not of generalized Roter type. We now discuss the results for particular value of the functions $f$, $h$ and $\psi$ step by step as follows:\\
\textbf{Step I:}  If $(h')^2-h h'' =0$, i.e., $h=c_1 e^{c_2 x^2}$, then $M$ is generalized Roter type  and also satisfies the $Ein(3)$ condition. We note that in this case fiber $\widetilde M$ is proper Roter type and thus $M$ is a proper generalized Roter type warped product manifold with proper Roter type fiber.\\
\textbf{Step II:} Again consider $-2 (f')^2 + f(2 f'' -1) =0$, i.e.,
$$f = \frac{e^{-\sqrt{c_1}(x^1+c_2)}\left(e^{\sqrt{c_1}(x^1+c_2)}+4 c_1\right)^2}{16 c_1^2} \ \mbox{ or} \ \ 
f = \frac{e^{-\sqrt{c_1}(x^1+c_2)}\left(1 + 4 c_1 e^{\sqrt{c_1}(x^1+c_2)}\right)^2}{16 c_1^2},$$
where $c_1$ and $c_2$ are arbitrary non-zero constants. Then the manifold satisfies the $Ein(2)$ condition and thus the manifold becomes proper Roter type. In this case fiber remains also Roter type. So $M$ is a warped product proper Roter type manifold with proper Roter type fiber.\\
\textbf{Step III:} Next consider $(\psi')^2 - 2 \psi \psi'' =0$, i.e., $\psi =\frac{(c_1 x^3 +2 c_2)^2}{4 c_2}$. Then $M$ is of constant curvature and in this case fiber is also of constant curvature.\\
\indent We now discuss a special case, when $f =(x^1)^2$, $h=c_2 Cos^2(x^2-2c_1)$ and $\psi = e^{x^3}$. Here the manifold $M$ is a special generalized Roter type and satisfies the $Ein(3)$ condition. In this case the fiber $\widetilde M$ is proper generalized Roter type and $Ein(3)$. Hence $M$ is a warped product proper generalized Roter type manifold with proper generalized Roter type fiber. \\
%%%%%%%%%%%%%%%%%%%%%%%%%%%%%%%%%%%%%%%%%%%%%%%%%%%%%%%%%%%%%%%%%%%%%%%%%%%%%%%%%%%%%%%%%%%%%%%%%%%%%%%%%%%%%%%%%
\noindent
\textbf{Conclusion:}
%=======================
The characterization of Roter type and generalized Roter type warped product manifolds is investigated along with their proper existence by suitable examples.\\
%%%%%%%%%%%%%%%%%%%%%%%%%%%%%%%%%%%%%%%%%%%%%%%%%%%%%%%%%%%%%%%%%%%%%%%%%%%%%%%%%%%%%%%%%%%%%%%%%%%%%%%%%%%%%%%%%
\noindent
\textbf{Acknowledgment:}
%=======================
The second named author gratefully acknowledges to CSIR, 
New Delhi (File No. 09/025 (0194)/2010-EMR-I) for the financial assistance. All the algebraic computations of Section \ref{exam} are performed with help of Wolfram Mathematica.
%%%%%%%%%%%%%%%%%%%%%%%%%%%%%%%%%%%%%%%%%%%%%%%%%%%%%%%%%%%%%%%%%%%%%%%%%%%%%%%%%%%%%%%%%%%%%%%%%%%%%%%%%%%%%%%%%

%%%%%%%%%%%%%%%%%%%%%

%%%%%%%%%%%%%%%%%%%%%%%%%%%%%%%%%%%%%%%%%%%%%%%%%%%%%%%%%%%%%%%%%%%%%%%%%%%%%%%%%%%%%%%%%%%%%%%%%%%%
\end{document}